\documentclass[leqno]{article}
\usepackage{latexsym,amsmath,amsfonts,mathrsfs,amssymb,amsthm}
\usepackage{yhmath}
\usepackage[usenames]{color}
\usepackage{graphicx}
\newcommand{\cupdot}{\mathbin{\mathaccent\cdot\cup}}

\def\R{\mathbb R}
\def\N{\mathbb N}

\def\S{\mathbb S}
\def\s{\sharp}

\def\e{\epsilon}
\def\d{\delta}

\def\be{\begin{equation}}
\def\ee{\end{equation}}
\def\bs{\backslash}
\def\qed{\hfill$\Box$\bigskip}
\def\nd{\noindent\textbf{Proof. }}

\def\lg{\langle}
\def\llg{\left\langle}
\def\rg{\rangle}
\def\rrg{\right\rangle}
\def\tb{\textcolor{blue}}

\numberwithin{equation}{section}
\newtheorem{lem}[equation]{Lemma}

\newtheorem{defn}[equation]{Definition}
\newtheorem{thm}[equation]{Theorem}
\newtheorem{cor}[equation]{Corollary}

\begin{document}
\bigskip

\centerline{\Large \textbf{On the dimension of $k$-medial axis for arbitrary closed set}}

\bigskip

\centerline{\large Xiangyu Liang}

\begin{center}School of mathematical sciences, Beihang University \\
maizeliang@gmail.com \end{center}

\vskip 1cm

\centerline {\large\textbf{Abstract}}

We prove that the $k$-medial axis of an arbitrary closed set in $\R^n$ is $n-k+1$-rectifiable (and hence of dimension at most $n-k+1$). This result gives a first stratification for medial axis of any closed set, which has been widely studied and used in pure and applied mathematics. This also answers a question proposed by Erd\"os \cite{Erd45}, and leads to more further interesting investigations (see the end of the article).

\bigskip


\bigskip


\section{Introduction}

In this article we give a first regularity result for $k$-medial axis of an arbitrary closed set in Euclidean spaces, and hence obtain a stratification for medial axis of any closed set. 

Given a closed set $E\subset\R^n$, denote by $M$ the set of points for which there is more than one nearest point in $E$. M is called the medial axis of the set $E$.

Due to the wide application in both theoretic and applied mathematics, medial axis and similar concepts such as central sets, skeletons, conflict sets, etc. have been extensively studied and used in various circumstances. Some classical literatures (far from exclusive) are the following: \cite{Hor83} in classical analysis, \cite{Erd45},  \cite{Frem97} and \cite{Pauc39} in geometric measure theory, \cite{Mil80}, \cite{Thom72} and \cite{Yom81} in differential geometry, \cite{FrPo80} in calculus of variations, and more recently, \cite{bis10} in conformal geometry. Besides, one can also find a tremondous number of literatures on application and research in applied mathematics, e.g. computational geometry, image analysis, etc. One well known example is proposed by Blum in 1967 \cite{Blum67} in shape analysis, the medial axis or the medial axis transform is a one dimensional graph extracted from a planar shape. More precisely, when the boundary of a domain is regular, it is known that the medial axis (or skeleton in bitmap image) can be stratified into the union of smooth curves and isolated singularities, and is a deformation retract of the region (cf. \cite{CCM97}). In this case, the shape of the region can be characterized by the geometry (type of lower dimensional singularities, stratification, etc.) of its medial axis, which is one dimensional less. This is used to compress the input data and expedite the extraction of image features. But if the region is not so regular, very little is known about the structure of the medial axis and to what point can the geometry of the medial axis reflect the shape of the region. 

In all the above subjects (except \cite{Erd45}), only (2-)medial axis of sets $E$ with good regularity (smooth, etc.) are investigated. In \cite{Pauc39}, Pauc proved that in the plane, $M_2$ is contained in the sum of countably many Jordan curves and $M_3$ is countable. In \cite{Erd45}, Erd\"os proved that for an arbitrary closed set $E$, its medial axis $M$ is of measure zero. Furthermore, he proposed to do some finer analysis on $k$-medial axis  $M_k$---points for which there is at least $k$ nearest points in generic position in $E$. (See Definition \tb{2.1}). Erd\"os proved that when $k=n+1$, $M_k$ is denumberable, and he conjectured that for general $1\le k\le n+1$, $M_k$ should be of dimension $n-k+1$.

In this article we give a complete answer to his conjecture: we show that for any closed set $E$, its $k$-medial axis should be $n-k+1$-rectifiable, hence is of dimension at most $n-k+1$. This gives a first regularity result for $k$-medial axis of an arbitrary closed set, and thus a stratification property for the well-studied medial axis. We also give examples of particular sets whose $k$-medial axes are of dimension strictly less than $n-k+1$, thus the conjecture is not completely true.

Based on the current result, we also propose some further possible related problem at the end of this article.


\textbf{Acknowledgement:} This work is partially supported by China's Recruitement Program of Global Experts, School of Mathematics and Systems Science, Beihang University, and National Natural Science Foundation of China (Grant No. 11871090).

\noindent\textbf{Some useful notation}

For any $x,y\in \R^n$, $d(x,y)$ is the Euclidean distance  between $x$ and $y$; for any set $E\subset \R^n$, $d(x,E)$ denotes the distance from $x$ to $E$: $d(x,E)=\inf\{d(x,y):y\in E\}$.

%
%
%
$B(x,r)$ is the open ball with radius $r$ and centered on $x$;

$\overline B(x,r)$ is the closed ball with radius $r$ and center $x$;

For three points $x,y,z\in \R^n$, $\angle(x,y,z)\in[0,\pi]$ denotes the smaller angle between the two vectors $x-y$ and $z-y$.

In this article, by dimension we always mean Hausdorff dimension.

%
%
%
%
\section{The main result}
Let $E$ be any closed set in $\R^n$. For any $x\in \R^n$, set $d(x)=d(x,E)$ the distance from $x$ to $E$. we call $y\in E$ a point closest to $x$ if $d(x,y)=d(x)$. 
Denote by $\phi(x)$ the set of all points closest to $x$. 

\begin{defn}For any positive integer $1\le k\le n+1$, denote by $M_k=M_k(E)$ the set of points for which $\phi(x)$ contains $k$ points not all in a $(k-2)$-dimensional affine plane -- that is, the $k$ points are in generic position. Hence $M_k\supset M_{k+1}$ for any $k$. $M_k$ is called the $k$-medial axis of $E$.
\end{defn}

Erd\"os conjectured in \cite{Erd45} that the Hausdorff dimension of $M_k$ is $n+1-k$. A particular case was proved in \cite{Erd45} : when $n=2$ and $k=3$, the set $M_3$ is countable, in particular its dimension is 0. 

Notice that for a general $k$, it is obvious that if $E$ is convex, then for $k\ge 2$, $M_k=\emptyset$, thus the dimension of $M_k$ is zero. Hence the conjecture is not true. However, we will prove the other direction, that is the dimension of $M_k$ is at most $n+1-k$:

\begin{thm}Let $E$ be a closed subset of $\R^n$. Then for $1\le k\le n+1$, the $k$-medial axis $M_k$ of $E$ is $n+k-1$-rectifiable. In particular, $M_k$ is of Hausdorff dimension at most $n+k-1$.
\end{thm}

Note that after this Theorem \tb{2.2}, we know that for every closed $E\subset \R^n$ and each $k\ge 2$, the $n$-dimensional Lebesgue measure of $M_k(E)$ is zero, and hence $M_1(E)=\R^n\bs(\cup_{k\ge 2}M_k(E))$ is of dimension $n$. Therefore Erd\"os' conjecture is true for the particular case $k=1$.

\smallskip

\noindent\textbf{Proof of Theorem \tb{2.2}.}

When $k=1$, this is trivial, since $M_k$ is a subset of $\R^n$, and hence its dimension is at most $n=n+1-k$. 

So now fix any integer $2\le k\le n+1$. 

Denote by $A$ the set of all subsets $a$ of the unit sphere $\S^{n-1}$ with $k$ elements that span a $k-1$-dimensional affine plane, that is, 
\be A=\{a\subset S^{n-1}: \s a=k\mbox{ and the elements of\ }a\mbox{ span a }k-1\mbox{-dimensional affine plane}\}.\ee
Then the $k$ points in $a$ span a $k-1$-dimensional affine plane. Define a distance $d_A$ on $A$: for any $a,b\in A$, $a=\{a_1,a_2,\cdots, a_k\}$ and $b=\{b_1,b_2,\cdots b_k\}$ with $a_i,b_j\in S^{n-1},1\le i,j\le k$, 
\be d_A(a,b)=\inf_{\sigma\in\Sigma_k}\sup_{1\le i\le k}d(a_i,b_{\sigma(i)}),\ee
where $\Sigma_k$ denotes the symmetric group of $\{1,2,\cdots k\}$.

For any $x\in M_k$, denote by $B(x)=\{\frac{y-x}{d(x)}:y\in \phi(x)\}$. Then $B(x)\subset S^{n-1}$. We say that $d_M(x,a)<\e$ if $B(x)$ has a subset $b\in A$ such that $d_A(b,a)<\e$.

Let $d,\e,\d>0$, $a\in A$. Let $M(a,d,\e,\d)$ denote the subset of $M_k$:
\be M(a,d,\e,\d)=\{x\in M_k: |d(x)-d|<\d,  d_M(x,a)<\e.\}\ee

The idea of the proof is the following : we would like to prove that for any $a,d$, the set $M(a,d,\e,\d)$ is $n+k-1$ rectifiable with locally finite measure for $\e,\d$ small (depending on $a$ and $d$). Then we decompose $M_k$ into a countable union of sets of the form $M(a,d,\e,\d)$, and we are done.

To prove the rectifiability of $M(a,d,\e,\d)$, we will use the basic fact (cf. \cite{Ma} Lemma 15.13) that, if $P$ is a $d$-dimensional plane in $\R^n$, $C$ is an open cone containing $P^\perp$, then if a set $E$ is such that for each $x\in E$, $E$ does not meet $x+C$, then $E$ is a Lipschitz graph over a setset of $P$.

\medskip

So fix any $a\in A$.

Denote by $P(a)$ the $k-1$ affine plane spanned by the $k$ elements $a_1,\cdots, a_k$ of $a$. Denote by $P$ the $k-1$ plane parallel to $P(a)$ which passes through the origin ,and $Q$ the $n+1-k$-plane which is orthogonal to $P$ and passes through the origin. For $r>0$, denote by $P_r$ the cone 
\be P_r=\left\{v\in \R^n:\frac{\mbox{dist}(v,P)}{||v||}<r\right\}.\ee

Denote by $\pi_P:\R^n\to P$ and $\pi_Q:\R^n\to Q$ the orthogonal projections from $\R^n$ to $P$ and $Q$ respectively. Then by definition of $P$ and $Q$, the vectors $\pi_Q(a_i),1\le i\le k$ are equal, and the lengths $||\pi_P(a_i)||,1\le i\le k$ are equal. Set $h=\pi_Q(a_1)(=\pi_Q(a_i),1\le i\le k)$, then $||h||<1$. For any non zero vector $w\in P$, set 
\be f(w)=\left[\sup_{1\le i\le k}\llg w,\pi_P(a_i)\rrg -\inf_{1\le i\le k}\llg w,\pi_P(a_i)\rrg \right].\ee

\begin{lem} $1^\circ$ For any non zero vector $w\in P$, $f(w)>0$;

$2^\circ$ Set $c=\inf\{f(w):w\in P, \frac12\le ||w||\le 1\}$. Then $c>0$.
\end{lem}

\nd $1^\circ$ If the conclusion of $1^\circ$ does not hold, then there exists $w\in P\bs \{0\}$ so that $f(w)=0$. By definition of $f$, the vectors $\pi_P(a_i),1\le i\le k$ belong to the same $k-2$ affine subspace $\{u\in P:\llg u,w\rrg =\llg w,\pi_P(a_1)\rrg \}$ of $P$, and thus $a_i=h+\pi_P(a_i),1\le i\le k$ belong to a $k-2$ affine plane of $\R^n$, which contradicts our assumption that $a_i,1\le i\le k$ do not belong to a same $k-2$ affine plane. (See Figure \tb{1} for an idea in $\R^3$: when 3 vectors span the 2-plane $P$, then the projections of any non-zero vector to these 3 vectors cannot be the same.)

\centerline{\includegraphics[width=0.5\textwidth]{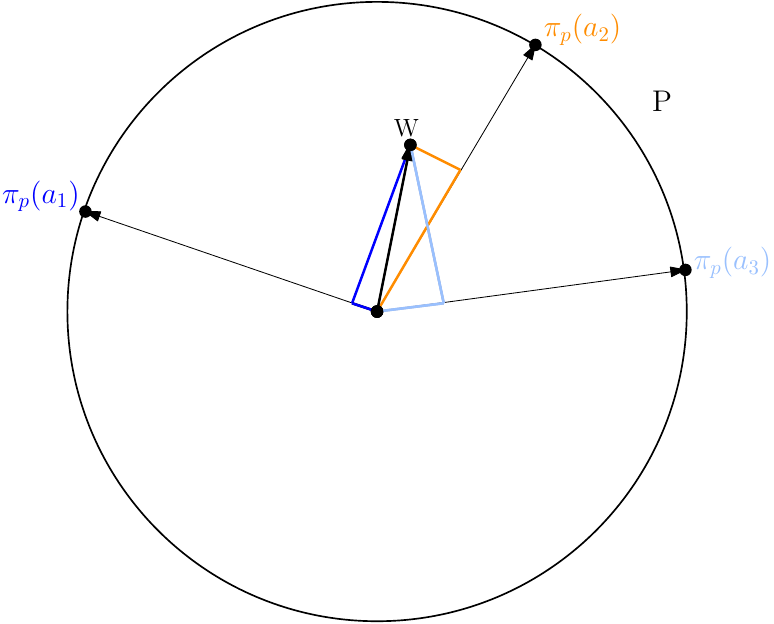}}
\nopagebreak[4]
\centerline{Figure 1}

$2^\circ$ It is easy to see that $f$ is continuous, and the set $\{w\in P,\frac12\le ||w||\le 1\}$ is compact. Thus the conclusion of $2^\circ$ follows directly from $1^\circ$.\qed

\begin{lem}For any $r<\frac 12$, and for any $v\in P_r\bs \{0\}$, we have
\be \sup_{1\le i\le k}\cos\angle(v,0,a_i)-\inf_{1\le i\le k}\cos\angle(v,0,a_i)>c.\ee
\end{lem}

\nd Take any $r<\frac 12$, and any  $v\in P_r\bs\{0\}$, we have
\be\cos\angle(v,0,a_i)=\frac{\llg v,a_i\rrg }{||v||}=\llg \frac{\pi_Q(v)}{||v||},h\rrg +\llg \frac{\pi_P(v)}{||v||},\pi_P(a_i)\rrg .\ee
Since $r<\frac 12$, $\frac{\pi_P(v)}{||v||}\in \{w\in P,\frac12\le ||w||\le 1\}$, hence by Lemma \tb{2.8 $2^\circ$}, 
\be\begin{split}&\sup_{1\le i\le k}\cos\angle(v,0,a_i)-\inf_{1\le i\le k}\cos\angle(v,0,a_i)\\
=&\sup_{1\le i\le k}\left\lg \frac{\pi_P(v)}{||v||},\pi_P(a_i)\right\rg -\inf_{1\le i\le k}\llg \frac{\pi_P(v)}{||v||},\pi_P(a_i)\rrg \\
>&c.\end{split}\ee

\qed

Now let us give a rough idea about what happens next (Figure \tb{2} below gives a geometric intuition): suppose that $0\in M_k$, with $d(0)=1$, $a\in A$ and $a\subset \phi(0)$. Then $0\in M(a,1,\e,\d)$ for any $\e>0$ and $\d>0$. We want to show that when $\e,\d$ are very small, then the set $M(a,1,\e,\d)$ does not intersect the cone $P_r\bs\{0\}$ in a small neighborhood $B(0,t)$ of 0.

Take any $v\in B(0,t)\cap P_r\bs \{0\}$. Lemma \tb{2.9} says that for there exists two vectors among the $a_i's$, , say $a_1$ and $a_2$, so that $\cos\angle(v,0,a_1)-\cos\angle(v,0,a_2)>c$. As a result, we know that $d^2(v)\le d^2(v, a_1)$, which is smaller than $d^2(v, a_2)-2||v||(\cos\angle(v,0,y_1)-\cos\angle(v,0,y_2))<d^2(v, a_2)-2c||v||$ by the cosine formula.

Then if $v\in M(a,d,\e,\d)$, this means that there exists $z\in B(a_1, \e)$ so that $v+d(v)z\in E$. But a simple geometry shows that $v+d(v)z\in B(0, 1-\frac 14c||v||)$, which is impossible because $d(0)=1$. Thus $M(a,1,\e,\d)\cap B(0,t)\cap P_r\bs\{0\}=\emptyset$.

\centerline{\includegraphics[width=0.5\textwidth]{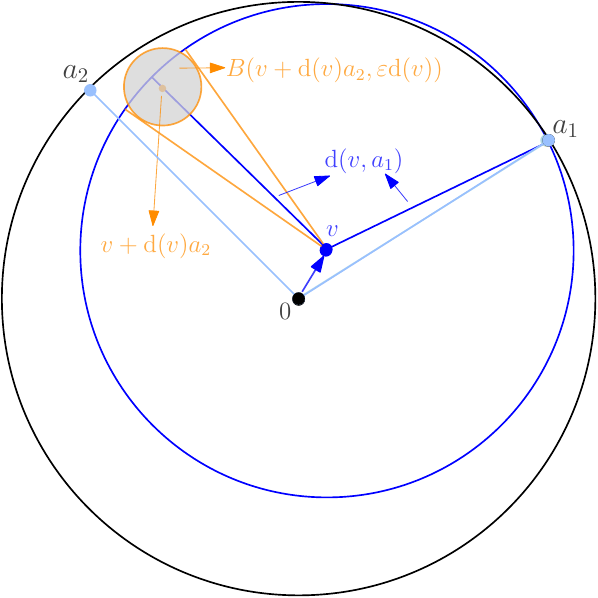}}
\nopagebreak[4]
\centerline{Figure 2}

Now we write down the above idea in detail in the following lemma, replacing 0 by a general point $x\in M_k$ so that $B(x)$ is very near $a$, and 1 by a general distance $d$. 

So fix any $d$.

\begin{lem}We can take $\e,\d, t,r>0$ small, which depends only on $a$ and $d$ such that for any $x\in M(a,d,\e,\d)$, $((P_r\bs\{0\})+x)\cap B(x,t)\cap M(a,d,\e,\d)=\emptyset$. 
\end{lem}

\nd Now let $\e,\d,r$ and $t$ be small, to be decided later. Take $x\in M(a,d,\e,\d)$, and $v\in P_r\cap B(0,t)$ which is non zero. By Lemma \tb{2.9}, without loss of generality, suppose that $\cos\angle(v,0,a_1)-\cos\angle(v,0,a_2)>c$. Let $y_1,\cdots, y_k\in B(x)$ be $k$ points such that 
 \be d(y_i,a_i)<\e.\ee
 
Now look at the point $x+v$.  Since $y_1\in B(x)$, the point $p=x+d(x)y_1\in E$. Hence 
\be d(x+v)=d(x+v,E)\le d(x+v,p).\ee Apply the cosine formula to the triangle with vertices $x,p$ and $x+v$, we have
\be d(x+v,p)^2=||v||^2+d(p,x)^2-2||v||d(p,x)\cos\angle(x+v,x,p).\ee
Notice that $\angle(x+v,x,p)=\angle(v,0,y_1)$, and $d(p,x)=d(x)$, therefore
\be d(x+v)^2\le d(x+v,p)^2=||v||^2+d(x)^2-2||v||d(x)\cos\angle(v,0,y_1).\ee

Now take any vector $z\in S^{n-1}$ such that $d(z,a_2)<\e$. Let us calculate the distance between $q=x+v+d(x+v)z$ and $x$. Apply again the cosine formula to the triangle with vertices $q,x$ and $x+v$:
\be d(x,q)^2=||v||^2+d(x+v,q)^2-2||v||d(x+v,q)\cos\angle(x,x+v,q).\ee

Notice that 
\be\angle(x,x+v,q)=\angle(-v,0,z)=\pi-\angle(v,0,z)\mbox{, and }d(x+v,q)=d(x+v),\ee 
hence
\be \begin{split}
d(x,q)^2=&||v||^2+d(x+v)^2-2||v||d(x+v)\cos(\pi-\angle(v,0,z))\\
=&||v||^2+d(x+v)^2+2||v||d(x+v)\cos\angle(v,0,z).\end{split}\ee

Combine \tb{(2.17) and (2.20)}, we have
\be \begin{split}
d(x,q)^2\le&2||v||^2+d(x)^2-2||v||d(x)\cos\angle(v,0,y_1)+2||v||d(x+v)\cos\angle(v,0,z)\\
=&2||v||^2+d(x)^2-2||v||d(x)\cos\angle(v,0,y_1)+2||v||d(x)\cos\angle(v,0,z)\\
&-2||v||d(x)\cos\angle(v,0,z)+2||v||d(x+v)\cos\angle(v,0,z)\\
=& 2||v||^2+d(x)^2+2||v||d(x)[\cos\angle(v,0,z)-\cos\angle(v,0,y_1)]\\
&+2||v||[d(x+v)-d(x)]\cos\angle(v,0,z)\\
\le&2||v||^2+d(x)^2+2||v||d(x)[\cos\angle(v,0,z)-\cos\angle(v,0,y_1)]+2||v|||d(x+v)-d(x)|\\
\le&d(x)^2+2||v||d(x)[\cos\angle(v,0,z)-\cos\angle(v,0,y_1)]+4||v||^2.
\end{split}\ee

Recall that $d(z,a_2)<\e$ and $d(y_1,a_1)<\e$. Hence when $\e$ is small enough (which depend only on $c$, and thus on $a$), we have
\be \cos\angle(v,0,z)-\cos\angle(v,0,y_1)>\frac c2,\ee
which yields
\be d(x,q)^2<d(x)^2-c||v||d(x)+4||v||^2.\ee

But $x\in M(a,d,\e,\d)$, which gives $d(x)>d-\d$. Therefore when $\d$ is small, which depends on $d$, $d(x)$ is no less than $\frac 12d$.  Also, since $0\ne v\in B(0,t)$, when $t$ is small enough, we have
\be d(x,q)^2\le d(x)^2-\frac12 cd||v||<d(x)^2.\ee
This means, for any $z\in S^{n-1}$ such that $d(z,a_2)<\e$, the point $q=x+v+d(x+v)z$ does not belong to $E$, because $d(x,q)<d(x,E)$. In other words, the set $B(x+v)$ is $\e$ far from the point $a_2$. In particular, $d_M(x+v,a)>\e$. Hence $x+v\not\in M(a,d,\e,\d)$. \qed

\begin{cor}
Fix any $a$ and $d$. Then there exist $\e=\e_{a,d}$ and $\d=\d_{a,d}$ (depending on $a$ and $d$) such that for any $\e\le\e_{a,d},\d\le\d_{a,d}$, the set $M(a,d,\e,\d)$ is $n+k-1$ rectifiable with locally finite measure.
\end{cor}

\nd This follows directly from Lemma \tb{2.13} and \cite{Ma} Lemma 15.13.\qed

Let us finish the proof of Theorem \tb{2.2}. 

For each $x\in M_k$, by definition, we can fix $k$ points $a_1^x, \cdots, a_k^x\in \phi(x)$ which do not lie in a $k-2$-dimensional affine subspace, such that $d(x,a_i^x)=d(x), \forall 1\le i\le k$. Let $a(x)=\{\frac{a_i^x}{d(x)},1\le i\le k\}\in A$.

Then we have
\be  M_k=\bigcup_{(a,d)\in A\times (0,\infty)}\{x: a(x)=a, d(x)=d\}.\ee

By Corollary \tb{2.25}, for each pair $(a,d)\in A\times (0,\infty)$, there exists $\e_{a,d}>0,\d_{a,d}>0$ such that the set $M(a,d,\e,\d)$ is $n-k+1$ rectifiable with locally finite measure. Denote by $Q_{a,d}=\{(a',d')\in A\times (0,\infty): |d'-d|<\d,d_M(a,a')<\e\}$. Then $Q_{a,d}$ is an open subset of $A\times (0,\infty)$, and the family $\{Q_{a,d},(a,d)\in A\times (0,\infty)\}$ forms an open cover of $A\times (0,\infty)$. Since the space $A\times (0,\infty)\subset \R^n\times \R$ is second countable, there exists a countable subcover $\{Q_{a_j,d_j}:j\in \N\}$. Then by \tb{(2.26)},
\be \begin{split}M_k&=\cup_{(a,d)\in A\times (0,\infty)}\{x: a(x)=a, d(x)=d\}\\
&=\cup_{j\in \N}\{x: (a(x),d(x)\in Q_{a_j,d_j}\}\\
&=\cup_{j\in \N}M(a_j,d_j,\e_{a_j,d_j},\d_{a_j,d_j}),
\end{split}\ee
Hence $M_k$ is a countable union of $n-k+1$ rectifiable sets, and thus is $n-k+1$-rectifiable.\qed

As a direct corollary, here is a ''stratification'' for medial axis for an arbitrary closed set in $\R^n$:

\begin{cor}Let $E\subset \R^n$ be closed, and for $1\le k\le n+1$, denote by $M_k$ the $k$-medial axis of $E$. Then for each $1\le k\le n$, we can decompose the $k$-medial axis $M$ into a disjoint union of sets
\be M_k=\cupdot_{i=0}^{n-k+1} L_i,\ee
where each $L_i=M_{n-i+1}\bs M_{n-i+2}$ is $i$-rectifiable.
\end{cor}

\section{Further questions}

Based on Theorem \tb{2.2}, it would be interesting to know which sets have $k$-medial axis with full dimension and which do not (denote by degenerate $k$-medial axis). To some extent, we believe that degenerate $k$-medial axis come from good geometric property (for example, constant curvature, convexity, etc.) of the set.

A particular case of this has been thoroughly studied in more general circumstances: Chebyshev sets--sets to which every point admit only one nearest point. Chebyshev sets are exactly convex sets in many classical cases, such as Eudlidean spaces, many Hilbert spaces and smooth Banach spaces. But the convexity of Chebyshev sets in general cases remains open. 

Back to degenerate $k$-medial axis, there are two kinds of them: 

$1^\circ$ sets whose $k$-medial axis is empty--we may call them $k$-Chebyshev sets; 

$2^\circ$ sets whose $k$-medial axis exists but is not of full dimension.

For the first case, when $k=1$, this is exactly the case for Chebyshev sets, and in Euclidean spaces Chebyshev sets are exactly convex sets. However when $k>1$, things might be more complicated. When $k=2$, a most intuitive kind of of $2-$Chebyshev sets are unions of two convex sets. But obviously this is not the only kind of $2$-Chebyshev sets. For example, one can imagine the union of two sets $E$ and $F$ of positive reach (see \cite{Fe59} for definition, classical examples are compact sets with smooth boundary), where $E^C$ lies in the reach of $F$, and $F^C$ lies in the reach of $E$.  However at the end of the day, we still expect $k$-Chebyshev sets admit simple geometries.

For the second case, we expect to have some curvature regularity for the original set $E$.

\renewcommand\refname{References}
\bibliographystyle{plain}
\bibliography{reference}

\end{document}